\documentclass[12pt, a4paper]{article}
\usepackage[monochrome]{color}
\usepackage{amsmath}
\usepackage{amsfonts}
\usepackage{amssymb}

\usepackage{epsfig}
\usepackage{graphicx}

\newcommand{\IR}{\ensuremath{\mathbb{R}}}

\newcommand{\IP}{\ensuremath{\mathbb{P}}}

\newcommand{\FF}{\ensuremath{{\cal{F}}}}

\newcommand{\IA}{\ensuremath{{\cal{A}}}}
\newcommand{\HH}{\ensuremath{\mathcal{H}}}

\newcommand{\BB}{\ensuremath{\mathcal{B}}}

\newcommand{\CC}{\ensuremath{\mathcal{C}}}

\newcommand{\C}{\mathcal C}

\newtheorem{theorem}{Theorem}[section]
\newtheorem{lemma}[theorem]{Lemma}
\newtheorem{hypothesis}[theorem]{Hypothesis}
\newtheorem{corollary}[theorem]{Corollary}

\newtheorem{proposition}[theorem]{Proposition}

\author{Erik Broman\footnote{ Uppsala Universitet, Sweden, {\tt broman@math.uu.se}.}, Federico Camia\footnote{VU University Amsterdam, The Netherlands, {\tt fede@few.vu.nl}. Research supported by NWO VIDI grant 639.032.916.}, Matthijs Joosten\footnote{VU University Amsterdam, The Netherlands, {\tt mjoosten@few.vu.nl}. Research supported by NWO grant 613.000.601.},
Ronald Meester\footnote{VU University Amsterdam, The Netherlands, {\tt rmeester@few.vu.nl}.}}
\title{Dimension (in)equalities and H\"older continuous curves in fractal percolation}

\begin{document}

\maketitle

\begin{abstract}
We relate various concepts of fractal dimension of the limiting set
$\mathcal C$ in fractal percolation to the dimensions of the set consisting of connected components larger than one point and its complement
in $\mathcal C$ (the ``dust''). In two dimensions, we also show that the set consisting
of connected components larger than one point is a.s.\ the union of non-trivial H\"older continuous
curves, all with the same exponent. Finally, we give a short proof of the fact that in two
dimensions, any curve in the limiting set must have Hausdorff dimension strictly larger than 1.
\end{abstract}



\section{Introduction and main results}\label{sec-intro}

In this paper we are concerned with a percolation model, first
introduced in~\cite{Mandelbrot}, which is known as Mandelbrot's
fractal percolation process and which can be informally described as follows.
For any integers $d \geq 2$ and $N \geq 2,$ we start by dividing the
unit cube $[0,1]^d \subset {\mathbb R}^d$ into $N^d$ closed subcubes
of equal size $1/N \times 1/N \times  \cdots \times 1/N$. Given $p \in [0,1]$ and a subcube, we retain the
subcube with probability $p$ and discard it with probability $1-p$.
This is done independently for every subcube of the partition.
Sometimes we adopt the terminology calling retained cubes {\em
black} and deleted cubes {\em white}. We define the random set
${\mathcal C}_N^1={\mathcal C}_N^1(d,p) \subset [0,1]^d$ as the
union of all retained subcubes.
Next consider any retained (assuming that ${\mathcal C}_N^1 \neq \emptyset$) subcube
$B$ in ${\mathcal C}_N^1.$ We repeat the described procedure on a smaller
scale by dividing $B$ into $N^d$ further subcubes, discarding or
retaining them as above. We do this for every retained subcube of
${\mathcal C}_N^1.$ This yields a new random set ${\mathcal C}_N^2 \subset {\mathcal C}_N^1.$
Iterating the procedure on every smaller scale yields an infinite
sequence of random sets $[0,1]^d \supset {\mathcal C}_N^{1} \supset {\mathcal C}_N^{2}
\supset \cdots $ and we define the limiting set
$$
{\mathcal C}_N:=\bigcap_{n=1}^{\infty} {\mathcal C}_N^n.
$$
We will hereafter suppress the $N$ in our notation and simply
write ${\mathcal C}$ for ${\mathcal C}_N$.

We will need a more formal definition of the model as well. Let
\[
 I^{n}_{k} := \left[\frac{(k-1)}{N^{n}},\frac{k}{N^{n}}\right],
\]
where $n\geq 1$ and $1 \leq k \leq N^{n}$.
For ${\bf k}=(k_1,\ldots,k_d),$ consider the subcube $D_{\bf k}^n$ of $[0,1]^d$ defined by
$D_{\bf k}^n:=I^{n}_{k_{1}} \times I^{n}_{k_{2}} \times \ldots \times I^{n}_{k_{d}}$, and let
$D:=\{D_{\bf k}^n: n\geq 1, 1 \leq k_{l} \leq N^{n} \}$. A cube $D_{\bf k}^n$ will sometimes be called a {\em level}-$n$
cube.
We define the sample space by
\[
 \Omega := \{0,1\}^{D},
\]
and denote an element of $\Omega$ by $\omega$. We let $\BB$ be the
Borel $\sigma$-algebra on $\Omega$ generated by the cylinders and
let $\IP_p$ denote the product measure on $\BB$ with density
$p\in[0,1]$, that is, we let $\IP_p(\omega(D_{\bf k}^n)=1)=p$
independently for every $D_{\bf k}^n\in D$. The limiting set is then
defined to be the intersection of all $D_{\bf k}^n \in D$ such that
$\omega(D_{\bf k}^n)=1$.

Let $CR([0,1]^d)$ denote the event that $\C$ contains a connected
component which intersects the left hand side $\{0\} \times
[0,1]^{d-1}$ of the unit cube and also intersects the right hand
side $\{1\} \times [0,1]^{d-1}$. In this case we say that a
\emph{left-right crossing} of the unit cube occurs.

We define the \emph{percolation function} $\theta_{N,d}$ by
\begin{equation}\label{crit-sheet-perc}
 \theta_{N,d}(p):=\IP_{p}(CR([0,1]^d)).
\end{equation}
The \emph{critical value} is defined as
\[
 \tilde{p}_c= \tilde{p}_{c}(N,d) := \inf \{p: \theta_{N,d}(p)>0\}.
\]
It has been shown in \cite{ChaChaDur88} that the phase transition in
Mandelbrot fractal percolation is non-trivial, i.e. $0 <
\tilde{p}_{c}(N,d)<1$. Furthermore it was discovered in
\cite{ChaChaDur88} that for $d=2$, $\theta_{N,d}(p)$ is
discontinuous at $\tilde{p}_c$ (see \cite{DekMee90} for an easy
proof). This was generalised in \cite{BroCam08} to all $d\geq 3$ and
$N$ large enough but the result is conjectured to hold for all $N$
in any dimension. (At this point we remark that as a corollary to
the proof of Theorem \ref{thm-main3} below, we obtain an explicit
bound for the size of the discontinuity at the critical
value $p_c$, defined below and conjectured to coincide with
$\tilde{p}_c$, in terms of the Hausdorff dimension of the set $\CC^c$,
also to be defined below.)

In $d=2$, the set ${\mathcal C}$ is a.s.\ totally disconnected for
$p<\tilde{p}_{c}(N,2)$. This is also known to be true in higher
dimensions for the same set of $N$ for which it is known that
$\theta_{N,d}(p)$ is discontinuous at $\tilde{p}_c$ (see
\cite{BC2}) and is conjectured to be true for all $d$ and
$N$.
It is therefore natural to work with the following critical
value:
\[
p_c(N,d):=\sup\{p: \CC \textrm{ is a.s.\ totally disconnected} \}.
\]
It is known (see \cite{BC2}) that for any $d\geq 2$ and $N\geq 2$,
$$\IP_p( {\mathcal C} \textrm{ is not totally disconnected})>0$$
if $p=p_c(N,d)$. Given this, it is an easy exercise to show that for $p \geq p_c(N,d)$, $\C \neq \emptyset$ implies that
the set $\C^c$ consisting of the union of all connected components larger than one point
is a.s.\ not empty.

We now fix $N,d\geq 2$ and assume that $p_{c}(N,d) \leq p <1 $.
For any point $x \in \CC$, let $\CC_x \subset \CC$  be the set of points $y \in \CC$ that are connected to $x$ in $\CC$.
We call $\CC_x$ the \emph{connected component} of $x$. It is known (see \cite{Mee92}) that for
$p \geq p_{c}(N,d)$ there exist a.s.\ uncountably many $x \in \CC$ such that $\CC_{x} = \{ x \}$.
We partition $\CC$ into two sets, $\CC^d:=\{x\in \CC: \CC_x=\{x\}\}$ and the aforementioned $\CC^c:=\CC\setminus \CC^d$. (To understand the
notation: {\em d} is short for ``dust", and {\em c} is short for ``connected".)

Before we can state our results we will need some more definitions. The reader is referred to
\cite{Falconer} for a general overview
of the subject of fractal sets.

A countable collection $\{B_i\}_{i=1}^{\infty}$
of subsets of ${\mathbb R}^d$ with diameter at most $\epsilon$ is called
an $\epsilon$-cover of $F$ if $F \subset \cup_{i=1}^{\infty}B_i$.
Define the \emph{$s$-dimensional Hausdorff measure} of $F$ as follows:
\[
{\cal H}^s(F):=\lim_{\epsilon \rightarrow 0} \inf\left\{ \sum_{i=1}^\infty diam(B_i)^s: \{B_i\}_{i=1}^{\infty}
\textrm{ is an } \epsilon\textrm{-cover of } F\right\}.
\]
The \emph{Hausdorff dimension} $\dim_{\HH} (F)$ of $F$ is defined as
\begin{equation} \label{eqn1}
\dim_{\HH}(F):=\inf \{s: {\cal H}^s(F)=0\},
\end{equation}
which also turns out to be equal to $\sup \{s: {\cal H}^s(F)=\infty\}$.
The Hausdorff dimension of the limiting set in fractal percolation is a.s.\ given
by the following equation, whose proof can be found in \cite{ChaChaDur88} or \cite{Falconer}, Proposition 15.4:
\begin{equation} \label{eqn:haus-dim-retain-set}
\dim_{\HH}({\mathcal C})=\left\{ \begin{array}{ll}
                     d+\frac{\log p}{\log N}& \textrm{if } {\mathcal C}\neq \emptyset,  \\
                     0 & \textrm{otherwise}.
                    \end{array}
\right.
\end{equation}

There are many other concepts of dimensionality and we will in particular use the following.
For a bounded set $F\subset {\mathbb R}^d$ let $M_{\delta}(F)$ be
the minimal number of closed
cubes of side length $\delta$ that is needed to cover $F$.

The \emph{Lower Box counting dimension} of $F\subset {\mathbb R}^d$ is given by
\[
\underline{\dim}_B(F):= \liminf_{\delta \rightarrow 0}\frac{\log M_{\delta}(F)}{-\log \delta},
\]
while the \emph{Upper Box counting dimension} of $F\subset {\mathbb R}^d$ is given by
\[
\overline{\dim}_B(F):= \limsup_{\delta \rightarrow 0}\frac{\log M_{\delta}(F)}{-\log \delta}.
\]
If $\underline{\dim}_B(F) = \overline{\dim}_B(F)$ then the common value is denoted $\dim_B(F)$ and called
the \emph{Box counting dimension} of $F$.
It is known (see e.g.\ \cite{Falconer}) that
for any bounded set $F\subset {\mathbb R}^d$
\begin{equation} \label{eqn4}
\dim_{\HH}(F)\leq   \underline{\dim}_B(F) \leq \overline{\dim}_B(F).
\end{equation}

The next two theorems contain our dimension results for fractal percolation.

\begin{theorem}
\label{thm-main3}
For $p_c(N,d)\leq p < 1$, we a.s.\ have
\begin{equation}\label{eqn-main3.1}
\dim_B ({\mathcal C}^c)=\dim_B ({\mathcal C})=\dim_{\HH}({\mathcal C}).
\end{equation}
If ${\mathcal C} \neq \emptyset$ then a.s.
\begin{equation}\label{eqn-main3.2}
\dim_{\HH} ({\mathcal C}^c)<\dim_{\HH} ({\mathcal C}),
\end{equation}
from which it easily follows that a.s.
\begin{equation}\label{eqn-main3.3}
\dim_{\HH} ({\mathcal C}^d)=\dim_{\HH} ({\mathcal C}).
\end{equation}
\end{theorem}

Note that if $p<p_c(N,d)$, ${\mathcal C}^c = \emptyset$ a.s.\ and so equation (\ref{eqn-main3.2})
still holds as long as $\dim_{\HH}({\mathcal C})>0$.

\begin{theorem}
\label{prop1}
For every $p$ there exists $1 \leq \beta=\beta(p) \leq d$ such that 
$$
\IP_p(\CC^c =\emptyset \text{ or } \dim_{\HH}(\CC^c)=\beta)=1.
$$
\end{theorem}

For $\epsilon>0$, let $\CC^{c,\epsilon}$ be the union of the
connected components of diameter at least $\epsilon$. The
following result suggests that the ``small components" of $\CC^c$
are the ones which actually determine its Box counting dimension.

\begin{proposition}
\label{thm-main4}
If ${\mathcal C}^c \neq \emptyset$, then
\[
{\mathbb E}_p [\underline{\dim}_B ({\mathcal C}^{c,\epsilon})]\leq D \dim_B ({\mathcal C}^c),
\]
where $D<1$ is independent of $\epsilon$.
\end{proposition}

\medskip
When $p \geq \tilde{p}_c$, it is natural to ask about the nature of the left-right crossings of the unit cube.
For $d=2$, it was shown in \cite{Mee92} that $\CC$ contains at least
one continuous curve crossing the square as soon as a connected
component crossing the square exists. It was later established in
\cite{Cha94} (again for $d=2$) that any curve in $\CC$ must have Hausdorff dimension strictly larger
than 1.

In this paper, focusing again on the two-dimensional version of the model, we take the issue of the existence
of continuous curves in $\CC$ much further, using the sophisticated machinery of Aizenman and Burchard \cite{AizBur99}.
Their paper deals with scaling limits of systems of random curves,  but we will show how their results can be
useful in the context of fractal percolation as well. This is perhaps somewhat surprising, since
the scaling limits in \cite{AizBur99} deal with convergence in distribution, whereas in the fractal context, the fractal
limiting set is an a.s.\ limit. The key will be a very careful comparison between convergence in the
weak sense of curves in an appropriate topology, and convergence in the a.s.\ sense of compact sets in another topology. 
From such a comparison, one can obtain information about the compact sets that make up the a.s.\ limit of the fractal
construction.

\begin{figure}
\label{figuur1}
\begin{center}
\includegraphics[width=8cm]{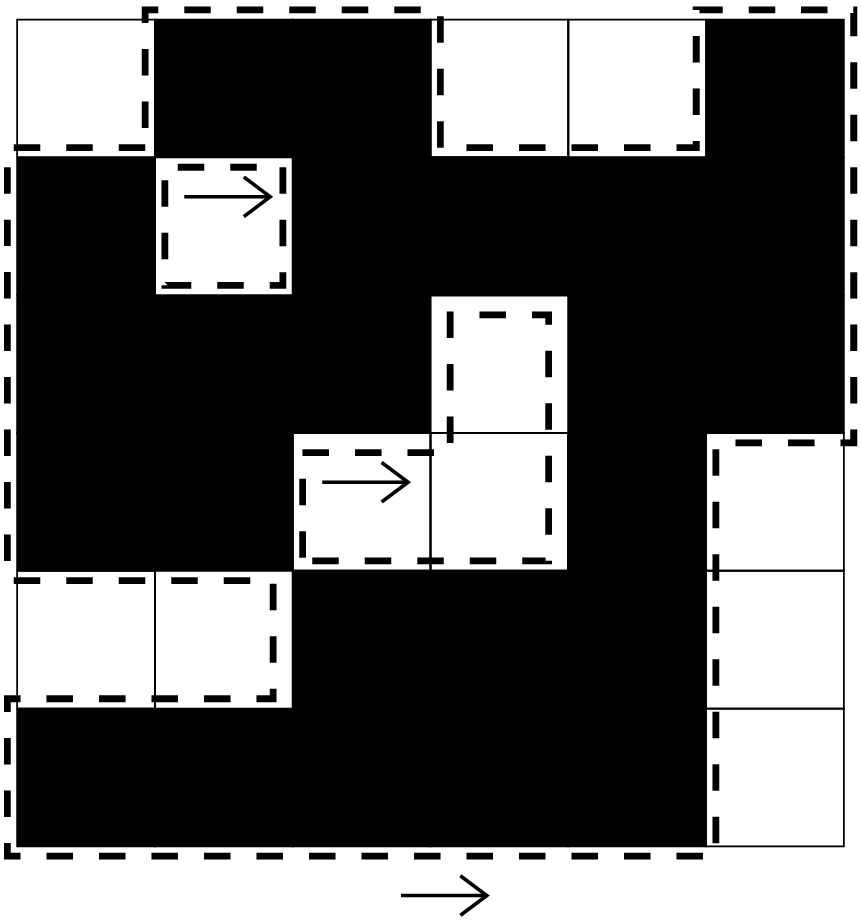}
\caption{The interface curves are drawn with broken lines. Arrows indicate the orientation.}
\end{center}
\end{figure}

In order to state our results, we need some definitions.
First of all, we define {\em interface curves} in the fractal
process. The complement $\mathbb{R}^2\setminus \CC^n$ consist of a
finite number of connected components, exactly one of which is
unbounded. The boundary of any such connected component can be split
into closed curves (loops). We call such loops {\em interface
curves} and denote by $\FF_n$ the collection of interface curves
after $n$ iterations of the fractal process. In order for our
interface curves to be uniquely defined, we orient them in such a
way that they have black (retained) squares on the left and white
(discarded) squares on the right, and assume that they turn to the
right at corners where two white and two black squares meet in a
checkerboard configuration, see Figure \ref{figuur1}. 

A connected
subset of an interface curve delimited by a starting and an ending
point will be called an \emph{interface segment}.

We continue with some general definitions concerning
curves, mostly taken from \cite{AizBur99} (see also \cite{CamNew06}).
We regard curves in $[0,1]^2$ as equivalence classes of continuous functions from $[0,1]$ to $[0,1]^2$
modulo strictly monotonic re-parametrizations. Below, $\gamma$ will
represent a particular curve and $\gamma(t)$ a particular parametrization
of $\gamma$. Denote by $\cal S$ the complete separable metric
space of curves in $[0,1]^2$ with metric
\begin{equation} \label{Distance}
\text{D} (\gamma_1,\gamma_2) := \inf
\sup_{t \in [0,1]} |\gamma_1(t)-\gamma_2(t)|,
\end{equation}
where the infimum is over all parametrizations
of $\gamma_1$ and $\gamma_2$.
The distance between two
sets ${\cal F}$ and ${\cal F}'$ of curves is defined by the Hausdorff metric induced by D, that is, 
$\text{Dist}({\cal F},{\cal F}') \leq \varepsilon$ if and only if
\begin{equation} \label{hausdorff-D}
\forall \, \gamma \in {\cal F}, \, \exists \,
\gamma' \in {\cal F}' \text{ with }
\text{D} (\gamma,\gamma') \leq \varepsilon
\text{ and vice versa}.
\end{equation}
The space $\Sigma$ of closed subsets of $\cal S$
with the metric Dist is also a complete
separable metric space.

The fractal process induces a probability measure $\mu_{n}$ on $\Sigma$, where $\mu_n$
denotes the distribution of $\FF_n$. With this notation, we can present our main result on continuous curves
in $\CC$.

\begin{theorem}
\label{thmnew1}
The sequence of measures $(\mu_n)$ has subsequential weak limits. Any such weak limit $\mu$
assigns probability 1 to curve
configurations in which all curves are H\"older continuous with the same exponent. The limiting set $\CC$ has the
same distribution as $g(\FF):=\cup_{\gamma \in \FF} {\rm Image}(\gamma)$, where $\FF$ is
a random set of curves distributed as $\mu$.
In other words, $\CC$ is distributed as the union of the images of the curves in a sample from a weak
limit of the $(\mu_n)$.
\end{theorem}

\medskip\noindent
{\bf Remark} Note that, although Theorem \ref{thmnew1} only claims the existence of subsequential weak limits
for $\mu_n$ as $n \to \infty$, Theorem \ref{thmnew1} combined with the existence of a unique a.s.\ limit for
$\CC^n$ implies that $g \circ \mu_n(\cdot) = \mu_n(g^{-1}(\cdot))$ has a unique weak limit.

\medskip
We now briefly discuss this result.
Since a single point is of course a H\"older continuous curve, the bare statement that a set is the union of H\"older
continuous curves is in itself close to being an empty statement. However, the curves in $\CC$ mentioned in
Theorem \ref{thmnew1} cannot be exclusively curves whose image is one point. One way to see this is to rephrase
the notion of weak convergence as follows (see also \cite{AizBur99}).
A sequence of probability measures $(\mu_n)$ on $\Sigma$ converges weakly to a probability $\mu$
measure on $\Sigma$ if and only if there exists a family of probability measures $\rho_n$ on
$\Sigma \times \Sigma$ such that the first marginal of $\rho_n$ is $\mu_n$, the second marginal
of $\rho_n$ is $\mu$ (for all $n$), and
$$
\int_{\Sigma \times \Sigma} {\rm Dist}(\FF_n, \FF)d\rho_n(\FF_n, \FF) \to 0
$$
as $n \to \infty$.
It is now also clear what happens to the points in the ``dust set"
$\CC^d$: these are accounted for as well in the theorem, since any
point $x \in \CC^d$ can be approximated by curves in $\FF_n$ whose
diameter and distance to $x$ converge to 0.

It is also possible to specialise to certain particular
curves. As an example, we discuss the {\em lowest crossing} in $\CC$ which we will
first properly define.
Condition on the existence in ${\CC}^n$ of a
left-right crossing of the unit square for all $n$, 
and consider the lowest interface segment $\sigma_n$ 
in ${\cal F}_n$ connecting the left and right side of the unit square. The closure of the region in the unit square
above $\sigma_n$ is a compact set (in the Euclidean topology) which decreases in $n$ and which therefore 
converges as $n \to \infty$.
The lowest crossing in $\CC$ is defined as the boundary of this limiting set.

\begin{theorem}
\label{thm:lowest-crossing} If $\CC$ contains a left-right
crossing of the unit square, then the lowest crossing in $\CC$ is a H\"older continuous curve.
\end{theorem}

The machinery of Aizenman and Burchard also allows for a quick proof, given in Section~\ref{sec-Hausdorff},
of the following result, first proved in \cite{Cha94}.

\begin{theorem}
\label{thmnew2} In two dimensions, there exists a constant
$\kappa>1$ such that all continuous curves in $\CC^c$ have Hausdorff
dimension at least $\kappa$.
\end{theorem}

\section{Proofs of Theorems \ref{thm-main3}, \ref{prop1} and Proposition \ref{thm-main4}}
\label{sec-dimension-bounds}

\medskip\noindent
\emph{Proof of Theorem \ref{thm-main3}.}
We start with the second equality of \eqref{eqn-main3.1}. Let $Z_n$ be the number of cubes of side length
$N^{-n}$ that are retained in ${\mathcal C}^n$. First observe that
$$
\frac{\log Z_n}{-\log N^{-n}} =
\frac{\log Z_n^{1/n}}{\log N}.
$$
It is well known from the theory of branching processes (see e.g. \cite{AthNey72,DekGri88}) that
$Z_n^{1/n}\rightarrow pN^d$ a.s.\ on the event ${\mathcal C}\neq \emptyset$. Therefore, for
a.e.\ $\omega$ such that ${\mathcal C}(\omega) \neq \emptyset$,
$$
\lim_{n \to \infty} \frac{\log Z_n}{-\log N^{-n}} = d + \frac{\log p}{\log N} ,
$$
and hence $\overline{\dim}_B(\CC) \leq d + \frac{\log p}{\log N}$.
Equations (\ref{eqn:haus-dim-retain-set}) and (\ref{eqn4}) then imply that
$\dim_B ({\mathcal C})=\dim_{\HH}({\mathcal C})$.

For the first equality of \eqref{eqn-main3.1}, let
$\{B_i\}_{i=1}^{M_{\delta}}$ be a cover of ${\mathcal C}^c$ using the minimal number $M_{\delta}$ of
closed cubes of side length $\delta$. For $A,B\subset \IR^d$, define $d(A,B):=\inf \{|x-y|:x\in A, y\in B\}$,
with $|\cdot|$ denoting Euclidean distance.
Assume that there exists $x\in {\mathcal C}$ such that $d(x,\bigcup_{i=1}^{M_{\delta}}B_i)>0$.
Then there must exist some $D_{\bf k}^n$ such that
$d(D_{\bf k}^n,\bigcup_{i=1}^{M_{\delta}}B_i)>0$ and $x\in D_{\bf k}^n$, which implies that $\omega(D^n_{\bf k})=1$ (i.e., $D^n_{\bf k}$ is retained).
However, because of the scale invariant construction of ${\mathcal C}$ and the fact that, for $p\geq p_c(N,d)$,
a.s.\ ${\mathcal C}$ is either empty or contains connected
components larger than one point (see \cite{BC2}), ${\mathcal C} \cap D_{\bf k}^n$ must
contain connected components larger than one point.
This contradicts the fact that $\{B_i\}_{i=1}^{M_{\delta}}$ is a cover of  ${\mathcal C}^c$ and shows that such an $x$ cannot exist.
Furthermore, since the union $\bigcup_{i=1}^{M_{\delta}}B_i$ is closed,
it follows that, if $d(x,\bigcup_{i=1}^{M_{\delta}} B_i)=0$ for $x \in \CC$, then $x \in \bigcup_{i=1}^{M_{\delta}}B_i$. Therefore, $M_{\delta}$
must be the minimal number of closed cubes of side length $\delta$ that covers ${\mathcal C}$.
This concludes the proof of \eqref{eqn-main3.1}.

Since ${\mathcal C}={\mathcal C}^c \cup {\mathcal C}^d$,
\eqref{eqn-main3.3} follows from \eqref{eqn-main3.2} and the fact that
(see \cite{Falconer})
\begin{equation} \label{eqn:max}
\dim_{\HH} ({\mathcal C})=\max(\dim_{\HH} ({\mathcal C}^c),\dim_{\HH} ({\mathcal C}^d)).
\end{equation}
We proceed therefore by proving \eqref{eqn-main3.2}, inspired by an argument suggested
by Lincoln Chayes to the second author.
Recall that ${\mathcal C}^{c,\epsilon}$ is the union of the connected components
of diameter at least $\epsilon$.
For $p \geq p_c(N,d)$, we have that
\begin{equation} \label{eqn5}
\dim_{\HH} ({\mathcal C}^c)=\sup_{\epsilon} \dim_{\HH} ({\mathcal C}^{c,\epsilon}),
\end{equation}
which is an easy consequence of the definition of Hausdorff dimension
(see \cite{Falconer}). Therefore, it suffices to find an upper bound of $\dim_{\HH} ({\mathcal C}^{c,\epsilon})$
which is uniform in $\epsilon$ and strictly smaller than $\dim_{\HH} ({\mathcal C})$.

We will now assume that $N\geq 5$ is odd. This assumption will make certain definitions easier to write down:
the reader can check that the squares $D_{\bf k}^n$ below would not be uniquely defined for even $N$.
We leave it to the reader to adapt the proof for all cases $N \geq 2$. 

For $D_{\bf k}^1$ such that
$(1/2,\ldots,1/2)\in D_{\bf k}^1$, let $B(D_{\bf k}^1;1):=[0,1]^d$
and $B(D_{\bf k}^1;3N^{-1})$ be the two cubes concentric to $D_{\bf
k}^1$ with side length $1$ and $3N^{-1}$ respectively. Let
$\varphi_{N,d}(p)$ be the probability that there exists a connected
component of $\C$ that crosses the ``shell'' $B(D_{\bf
k}^1;1)\setminus B(D_{\bf k}^1;3N^{-1})$. It follows from \cite{BC2}
that $\varphi_{N,d}(p)>0$ whenever $p \geq p_c(N,d)$.

Assuming that ${\mathcal C}^c \neq \emptyset$, fix
$\epsilon>0$ such that ${\mathcal C}^c$ contains at least one
component of diameter larger than $\epsilon$ and let $l$ be such
that $N^{-l+1} \leq \epsilon /d$. Consider a cube $D_{\bf k}^n$ for
$n \geq l$ which is intersected by a component of ${\mathcal C}^c$
of diameter larger than $\epsilon$. Let $B(D_{\bf k}^n;3N^{-n})$ and
$B(D_{\bf k}^n;N^{-n+1})$ be two cubes which are concentric to
$D_{\bf k}^n$ and have side length $3N^{-n}$ and $N^{-n+1}$
respectively. Obviously, for $D_{\bf k}^n$ to be intersected by a
connected component of diameter larger than $\epsilon$, there must
be a crossing of the shell $B(D_{\bf k}^n;N^{-n+1})\setminus
B(D_{\bf k}^n;3N^{-n})$. (Note that, depending on the
position of $D_{\bf k}^n$, it is possible that $B(D_{\bf k}^n;3N^{-n})$ and/or $B(D_{\bf
k}^n;N^{-n+1})$ are only partially contained in $[0,1]^d$.)

We will now construct a specific cover of $\CC^{c,\epsilon}$ which we will use
in our estimate for its Hausdorff dimension.
Let $W_n$ denote the set of cubes $D_{\bf k}^n$ with the following two properties:
\begin{itemize}
\item The intersection of all retained cubes of level $n$ and higher contains a crossing of the shell $B(D_{\bf k}^n;N^{-n+1})\setminus B(D_{\bf k}^n;3N^{-n})$.
In other words, if we would make all cubes black until level $n-1$
(inclusive), then there would be a connected component in $\CC$ crossing the shell $B(D_{\bf k}^n;N^{-n+1})\setminus B(D_{\bf k}^n;3N^{-n})$,
\item $D_{\bf k}^n$ is retained, that is, $\omega(D_{\bf k}^n)=1$.
\end{itemize}
By scale invariance and independence between the two conditions, we have that $\IP_p(D_{\bf k}^n\in W_n) \leq
p\varphi_{N,d}(p)$. The inequality is due to a boundary
effect since, as mentioned earlier, $B(D_{\bf k}^n;3N^{-n})$ and/or
$B(D_{\bf k}^n;N^{-n+1})$ need not be completely contained in
$[0,1]^d$.

\begin{figure}[t]
\label{figuur2}
\begin{center}
\includegraphics[width=8cm]{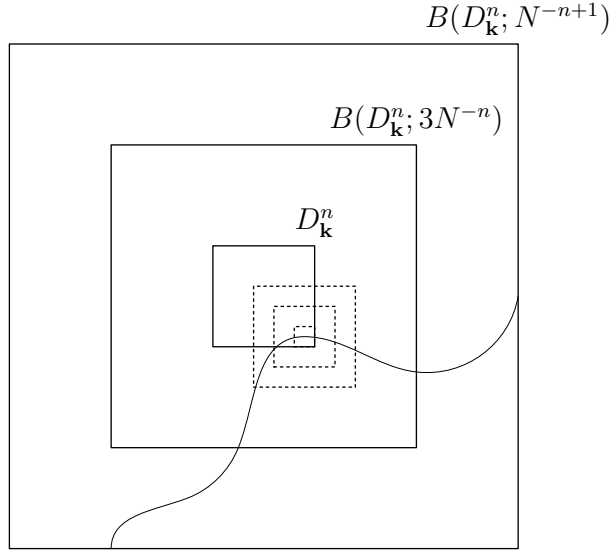}
\caption{The cube $D_k^{n}$ is in $V_{n}$ and $D_k^{n+1}$,
drawn with broken lines, belongs to $V_{n+1}$. Note that the corresponding shells are disjoint.}
\end{center}
\end{figure}

For a given cube $D_{\bf k}^n$, let $B^m$ denote the level-$m$ cube which contains $D_{\bf k}^n$, where $m \leq n$ (with
$B^n=D^n_{\bf k}$). 
We make two observations:
\begin{enumerate}
\item  If $D_{\bf k}^n$ has a non-empty intersection with $\CC^{c,\epsilon}$, then we have that $B^m \in W_m$, for all $m=l, l+1,\ldots, n$.
\item The events $\{B^m \in W_m\}$ form
a collection of independent events; see Figure \ref{figuur2}. 
\end{enumerate}
This motivates us to define $V_n$ as the collection of cubes $D_{\bf k}^n$ for which the corresponding cubes 
$B^m$ are in $W_m$, for all $m=l, l+1,\ldots, n$. From observation 1, we have that the collection $V_n$ forms a cover of $\CC^{c,\epsilon}$.
We can now write, using observation 2,
\begin{eqnarray*}
\IP_p(D_{\bf k}^n \in V_n) &=& \IP_p\left(\bigcap_{m=l}^n B^m \in W_m\right)\\
&=& \prod_{m=l}^n \IP_p(B^m \in W_m)\\
&\leq & (p\varphi_{N,d}(p))^{n-l+1}.
\end{eqnarray*}

Using Fatou's lemma and the fact that the collection of cubes in $V_n$ covers $\CC^{c,\epsilon}$, we obtain
(writing $|V_n|$ for the number of cubes in $V_n$)
\begin{eqnarray}
 {\mathbb E}_p\left({\cal H}^s({\mathcal C}^{c,\epsilon})\right) \nonumber
& \leq & \liminf_{n \to \infty} {\mathbb E}_p\left( \sum_{D_{\bf k}^n\in {V}_n} diam(D_{\bf k}^n)^s  \right) \nonumber \\
&=& \liminf_{n \to \infty} (\sqrt{d}N^{-n})^s{\mathbb E}_p(|V_n|) \nonumber \\
&  \leq & \liminf_{n \to \infty} d^{s/2}N^{-sn}(p\varphi_{N,d}(p)N^d)^{n-l+1} \nonumber \\
&=&  d^{s/2}N^{-d(l-1)}p\varphi_{N,d}(p) \lim_{n \to \infty} N^{n(d+\frac{\log(p\varphi_{N,d}(p))}{\log N}-s)}. \label{eqn-Hausdorff-limit}
\end{eqnarray}
The limit in \eqref{eqn-Hausdorff-limit} is finite if and only if
\[
s \geq d+\frac{\log(p\varphi_{N,d}(p))}{\log N},
\]
showing that
\[
\dim_{\HH}({\mathcal C}^{c,\epsilon}) \leq d+\frac{\log(p\varphi_{N,d}(p))}{\log N} \ \ \textrm{ a.s.}
\]
It follows from (\ref{eqn5}) that
\[
\dim_{\HH}({\mathcal C}^{c}) \leq d+\frac{\log(p\varphi_{N,d}(p))}{\log N} \ \ \textrm{ a.s.}
\]
and since $\varphi_{N,d}(p)<1$ the result follows from this and (\ref{eqn:haus-dim-retain-set}). \hfill $\Box$

\medskip\noindent
\emph{Proof of Theorem \ref{prop1}.} 
If $\IP_p(\C^c=\emptyset)=1$ then there is nothing to prove, so we assume that $\IP_p(\C^c \neq \emptyset)>0$.

Let $Z_n$ be the number of
retained cubes after $n$ steps of the fractal construction procedure
and let $B_1,\ldots,B_{Z_n}$ denote the retained cubes. If the event
$\{\dim_{\HH}({\mathcal C}^c)\geq \alpha\}$ occurs, then
$\{\dim_{\HH}({\mathcal C}^c\cap B_k)\geq \alpha\}$ for at least one $k=1,\ldots,Z_n$
(see, e.g., \cite{Falconer}). Therefore,
\[
\IP_p(\dim_{\HH}({\mathcal C}^c)\geq \alpha \Big{|}Z_n=l)
=1-\prod_{k=1}^l(1-\IP_p(\dim_{\HH}({\mathcal C}^c\cap B_k)\geq \alpha \Big{|} B_k \subset \C^n).
\]
However, because of scale invariance, $\dim_{\HH}({\mathcal
C}^c\cap B_k)$, conditioned on the event that $B_k \subset \C^n$,
must have the same distribution as $\dim_{\HH}({\mathcal C}^c)$ and
so in fact
\[
\IP_p(\dim_{\HH}({\mathcal C}^c)\geq \alpha \Big{|}Z_n=l)
=1-(1-\IP_p(\dim_{\HH}({\mathcal C}^c)\geq \alpha))^l.
\]
We can now write
\begin{eqnarray*}
\IP_p(\dim_{\HH}({\mathcal C}^c)\geq \alpha)
&  = & \sum_{l=1}^{N^{dn}}\IP_p(\dim_{\HH}({\mathcal C}^c)\geq \alpha \Big{|} Z_n=l)\IP_p(Z_n=l) \\
&  =& \sum_{l=1}^{N^{dn}}\left(1-(1-\IP_p(\dim_{\HH}({\mathcal C}^c)\geq \alpha))^l\right)\IP_p(Z_n=l)  \\
&  \geq &[1-(1-\IP_p(\dim_{\HH}({\mathcal C}^c)\geq \alpha))^n]\IP_p(Z_n\geq n).
\end{eqnarray*}
This last quantity is bounded below by
\begin{equation}
[1-(1-\IP_p(\dim_{\HH}({\mathcal C}^c)\geq \alpha))^n]\IP_p(Z_n\geq n \Big{|} \CC \neq \emptyset) \IP_p(\CC \neq \emptyset). \label{eqn11}
\end{equation}

As mentioned above,  a.s.\ $Z_n^{1/n}\rightarrow pN^d> 1$ as $n \to
\infty$ if $\CC \neq \emptyset$. Therefore, if
$\IP_p(\dim_{\HH}({\mathcal C}^c)\geq \alpha)>0$, by taking the limit of \eqref{eqn11}
as $n \to \infty$, we conclude that $\IP_p(\dim_{\HH}(\CC^c) \geq \alpha) \geq \IP_p(\CC
\neq \emptyset)$. If $\alpha>0$, it follows that in
fact $\IP_p(\dim_{\HH}(\CC^c) \geq \alpha) = \IP_p(\CC \neq
\emptyset)$.
Letting $\phi(\alpha)= \IP_p(\dim_{\HH}(\CC^c) \geq \alpha)$, we have that $\phi(0)=1$
and either $\phi(\alpha)=0$ or $\phi(\alpha)=\IP_p(\CC \neq \emptyset)$ when $\alpha>0$.

Now observe that ${\mathcal C}^c \neq \emptyset$ implies
$\dim_{\HH}({\mathcal C}^c) \geq 1$ (see, e.g., Proposition~4.1 of \cite{Falconer}).
Since the inverse implication is obvious, we conclude that $\phi(1)=\IP_p(\C^c\neq \emptyset)>0$
(where the last inequality follows from the assumption made at the beginning of the proof) and
that $\IP_p(\C^c\neq \emptyset)=\IP_p(\C\neq \emptyset)$.
Hence, since $\{\C^c\neq \emptyset\} \subset \{\C\neq \emptyset\}$, we obtain that
\begin{equation}
\label{zelfde}
\{\C^c\neq \emptyset\}= \{\C\neq \emptyset\},
\end{equation} 
up to a set of probability 0.

Moreover, since $\phi(\alpha)=\IP_p(\CC \neq \emptyset)$ for at least $\alpha \leq 1$,
and $\{ \dim_{\HH}({\mathcal C}^c)\geq \alpha \} \subset \{ {\mathcal C} \neq \emptyset \}$
for $\alpha>0$, we also obtain that $\{ \dim_{\HH}({\mathcal C}^c)\geq \alpha \} = \{ {\mathcal C} \neq \emptyset \}$,
up to a set of probability 0, for every $\alpha>0$ such that $\IP_p(\dim_{\HH}({\mathcal C}^c) \geq \alpha)>0$.

Writing
$$
\beta:=\sup\{\alpha:\IP_p(\dim_{\HH}({\mathcal C}^c)\geq \alpha)>0 \},
$$
we conclude from the above observations that, up to a set of probability 0,
if ${\mathcal C}^c \neq \emptyset$ then $\dim_{\HH}({\mathcal C}^c) = \beta$,
with $\beta \geq 1$ as a consequence of the assumption that $\IP_p(\dim_{\HH}(\CC^c) \neq \emptyset) > 0$
and the fact that $\dim_{\HH}({\mathcal C}^c) \geq 1$ as soon as ${\mathcal C}^c \neq \emptyset$. \hfill $\Box$

\medskip\noindent
Theorem \ref{prop1} and the proof of Theorem \ref{thm-main3} have an interesting
corollary which links $\dim_{\HH}({\mathcal C}^c)$ to the discontinuity at the critical
point $p_c(N,d)$. (This is the corollary that was announced in the introduction.)

\begin{corollary} \label{cor1}
Let $\Delta$ denote the a.s. Hausdorff dimension of ${\mathcal C}^c$ when
$p=p_c(N,d)$ and ${\mathcal C}^c \neq \emptyset$. Then,
\[
\varphi_{N,d}(p_c(N,d))\geq \frac{1}{N^{d-\Delta}}.
\]
\end{corollary}

\medskip\noindent
{\em Proof.} Let $p_b(N,d):=\inf\{p\leq
1:\varphi_{N,d}(p)>0\}$. Theorem~4.1 of \cite{BC2} and the
observation preceding it show that $p_b(N,d)=p_c(N,d)$. Moreover, it
follows from \cite{BC2} that $\varphi_{N,d}(p_c(N,d))>0$.
Combining these observations with the last line of the proof of
Theorem \ref{thm-main3}, we obtain that, for $p\geq p_c(N,d)$,
\[
\varphi_{N,d}(p) \geq p\varphi_{N,d}(p) \geq \frac{1}{N^{d-\Delta}}.
\]
\hfill $\Box$

\medskip\noindent
{\bf Remark} For Theorem \ref{thm-main3} we use the result from \cite{BC2} that for $p=p_c(N,d)$, $\IP_p(\CC^c \neq \emptyset)>0$.
It is possible to prove the result without this prior knowledge, as follows. We can start with the observation from \cite{BC2}
that $p_c=p_b$. We can then prove Theorem \ref{thm-main3} in the case of $p>p_c$ and from the last line of that proof we get that
$\varphi_{N,d}(p) \geq \frac{1}{N^{d-1}}$, using Corollary \ref{cor1} and the fact that the Hausdorff-dimension of a connected set
consisting of more than one point is at least 1 (see, e.g., Proposition 4.1 of \cite{Falconer}). Using this uniform bound and a
right-continuity argument similar to the ones in \cite{BC2}, we conclude that in fact $\varphi_{N,d}(p_c) \geq \frac{1}{N^{d-1}}$.
Hence we can {\em conclude} that for $p=p_c(N,d)$, $P_p(\CC^c \neq \emptyset)>0$, and go through the proof once more to obtain
the same result as above. 

\medskip\noindent
\emph{Proof of Proposition \ref{thm-main4}.} We again use the fact
that ${\mathcal C}^{c,\epsilon}$ is contained in the union of the
cubes in $V_n$, defined in the proof of Theorem~\ref{thm-main3},
and therefore, $M_{N^{-n}}= M_{N^{-n}}({\mathcal C}^{c,\epsilon}) \leq |V_n|$.
First observe that
\[
\liminf_{\delta \rightarrow 0} \frac{\log M_{\delta}}{-\log \delta}
=\liminf_{n \rightarrow \infty} \frac{\log M_{N^{-n}}}{-\log N^{-n}}.
\]
By Fatou's lemma and Jensen's inequality we get along the same lines as in the last part of the proof of Theorem \ref{thm-main3}
(using the same $l$ as in that proof) that
\begin{eqnarray*}
\lefteqn{{\mathbb E}_p\left(\liminf_{n \rightarrow \infty} \frac{\log M_{N^{-n}}}{-\log N^{-n}}  \right)} \\
&  \leq& \liminf_{n \rightarrow \infty} {\mathbb E}_p\left(\frac{\log M_{N^{-n}}}{-\log N^{-n}}  \right)
\leq \liminf_{n \rightarrow \infty} \frac{\log {\mathbb E}_p (M_{N^{-n}})}{-\log N^{-n}}\\
&  \leq & \liminf_{n \rightarrow \infty} \frac{\log {\mathbb E}_p \left(|V_n|  \right)}{n\log N}
\leq \liminf_{n \rightarrow \infty} \frac{\log (p\varphi_{N,d}(p)N^d)^{n-l+1}}{n\log N}\\
&  =& \frac{\log (p\varphi_{N,d}(p)N^d)}{\log N}
=d+\frac{\log(p\varphi_N(p))}{\log N}.
\end{eqnarray*}
Since $\varphi_{N,d}(p) < 1$ it follows that $d+\log(p \varphi_{N,d}(p))/\log N < d + \log p /\log N$. \hfill $\Box$

\section{Proof of Theorems \ref{thmnew1} and \ref{thm:lowest-crossing}}
In this section and the next one we work in two dimensions, that is, $d=2$. We will be using the machinery in
\cite{AizBur99} concerning scaling limits of systems of random curves.

One part of the argument is to apply some of the machinery developed in \cite{AizBur99} in order to show that the sequence of measures $\mu_{n}$ defined in the introduction has subsequential weak limits, and we first deal with this issue. After that, we combine the existence of weak limits with the a.s.\ limit behaviour of the fractal process in order to draw the final
conclusions.

Let $B(x,r)$ denote a closed square centered at $x$ with side length $r$ and $B^{\circ}(x;r)$ its
interior. For $R>r$, let $A(x;r,R):=B(x,R)\setminus B^{\circ}(x,r)$ be an annulus. The basic estimate is the following, from which everything else will follow.

\begin{lemma}\label{lem-disjoint-crossings} Let $p \geq p_c(N,2)$.
 There exists a sequence $\lambda(1), \lambda(2),\ldots$ with
$\lim_{k \to \infty} \lambda(k)=\infty$ and finite constants $K_k$
such that the following bound holds uniformly
for all $r \leq R \leq 1$ with $r$ small enough, and all $x$:
\begin{equation}\label{bnd-disjoint-crossings}
\limsup_{n\to\infty} \IP_p(\FF_n\text{ contains } k \text{ disjoint
crossings of } A(x;r,R)) \leq K_k
\left(\frac{r}{R}\right)^{\lambda(k)}.
\end{equation}
\end{lemma}

\medskip\noindent
{\em Proof.} We are looking for a collection of mutually disjoint
annuli (all contained in $B(x,R)$ and ``surrounding" $B(x,r)$) of
the form $A(y;2N^{-n}, 4N^{-n})$, where $y$ is a corner point of
some square $D^n_{\bf k}$ of the fractal construction. It is not
hard to see that for any $x \in [0,1]^2$, we can find such a
collection with at least $M:=c\log(R/r)$ elements, for a suitable
[uniform ???] positive constant $c$. For small enough $r$, $M$ is at least 1.
For a given $x$, we denote these
annuli by $A_{n_1}, A_{n_2},\ldots, A_{n_M}$, where the indices $n_1
< n_2 < \cdots < n_M$ refer to the $n$ associated with the annuli.
The idea is that if there are $k$ disjoint crossings of
$A(x;r,R)$, then each annulUS $A_{n_i}$ must
also be crossed by $k$ disjoint crossings in $\FF_n$. This is
exponentially unlikely in the number of annuli, as we will show now.

\begin{figure}[t]
\label{figuurtje}
\begin{center}
\includegraphics[width=8cm]{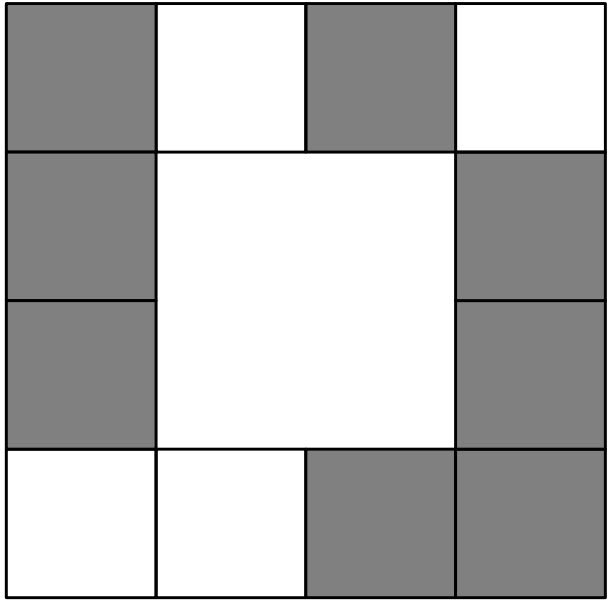}
\caption{The annulus $A_{n_1}$ with 3 level-$n_1$ components, drawn in dark grey.}
\end{center}
\end{figure}

We first consider the annulus $A_{n_1}$. Let $n > n_1$, and perform the fractal process until level $n_1$ (inclusive).
The annulus $A_{n_1}$ consists of 12 level-$n_1$ squares (Figure \ref{figuurtje}), some of which are in $\CC^{n_1}$
and some of which are not. Now observe that the collection of black (retained) squares of level $n_1$ in $A_{n_1}$
is partitioned into at most 6 ``level-$n_1$ components", where two neighbouring retained level-$n_1$ squares are
in the same component if they share an edge: see Figure \ref{figuurtje} for an illustration of the event that the
annulus contains 3 such level-$n_1$ components.

Now let $k$ be large enough (exactly how large will become clear soon).
If $A_{n_1}$ is, after $n$ iterations, crossed by $k$ disjoint interface segments, then at least
one of the level-$n_1$ components must be crossed by at least $k/6$ such segments, since there
are at most 6 components to accomodate these crossings (we should write $\lfloor k/6 \rfloor$ but
we ignore these details for the sake of notational convenience). Between the interface crossings,
one has alternating black (retained) and white (discarded) crossings of the annulus. This implies
that the component with at least $k/6$ interface crossings is, after $n$ iterations, also crossed
by at least $k/12-1$ white crossings which are ``between" the ``first" and ``last" of the $k/6$
interface crossings of said component.

The point of considering the level-$n_1$ components introduced above is that they are disjoint
and separated by vacant squares, so that the fractal constructions inside them from level $n_1$
on are independent of each other, and none of the interface segments crossing $A_{n_1}$ can
intersect more than one of them. Therefore, adding extra retained level-$n_1$ squares in $A_{n_1}$
does not affect the $k/12-1$ white crossings in the component with at least $k/6$ interface crossings.
This implies that the probability of having at least $k$ interface crossings in $A_{n_1}$ after $n$ iterations
is bounded above by the probability of having at least $k/12-1$ white crossings after $n$ iterations,
conditioned on having full retention in $A_{n_1}$ up to level $n_1$.

However, scaling tells us that, if we condition on retention until level $n_1$, the probability
in question is the same as the probability of having at least $k/12-1$ disjoint white components
crossing $\bar{A} := A((0,0);2,4)$ when we perform $n-n_1$ iterations of the fractal process in
$[-2,2]^2$ rather than in $[0,1]^2$, seen as the union of 16 independent fractal processes on the
16 unit squares making up $[-2,2]^2$. For these white crossing components we can use the BK inequality
(see \cite{Mee92} - a similar BK inequality for black crossing is not available) and deduce that the
probability of having $k/12-1$ of such white crossing components is bounded above by the probability
of having at least one, raised to the power $k/12-1$. This then finally leads to the estimate that
$$
\IP_p(\FF_n \text{ contains } k \text{ disjoint crossings of } A(x;r,R))
$$
is bounded above by
\begin{equation}
\label{mooi}
\IP_p(\bar{A} \text{ is crossed by a white component after } n-n_1 \text{ steps})^{k/12-1}.
\end{equation}

If there is a white component crossing $\bar{A}$, then there is no black circuit surrounding the origin in $\bar{A}$.
The probability of having such a black circuit after $n-n_1$ iterations is at least as large as the probability to
have such a circuit in the limit, and by the weak RSW theorem for fractal percolation in \cite{DekMee90} and the FKG
inequality, we have that for $p \geq p_c(N,2)$ this probability is strictly positive. It follows that there exists
$\alpha <1$ such that (\ref{mooi}) is bounded above by $\alpha^{k/12-2}$, uniformly in $n > n_1$.

Next we consider $A_{n_1}$ and $A_{n_2}$ simultaneously. Take $n > n_2$.
The probability to have $k$ interface crossings in $A_{n_1}$ and also in $A_{n_2}$ is the
probability that this happens in $A_{n_1}$ multiplied by the probability that this happens in $A_{n_2}$ conditioned on the fact that
it happens in $A_{n_1}$. We can treat this conditional probability exactly as above: we can change
the conditioning into one involving complete retention inside $A_{n_2}$ until level $n_2$,
to get an upper bound. It follows that the probability that both $A_{n_1}$ and $A_{n_2}$
have $k$ interface crossings is bounded above by the square of the individual bounds, that is, by $(\alpha^{k/12-1})^2$.

We continue in the obvious way now, leading to the conclusion that for $n > n_M$, the probability that
{\em all} annuli $A_{n_1}, \ldots, A_{n_M}$ are crossed by $k$ interface crossings, is bounded above by
$$
\left(\alpha^{k/12-1}\right)^{c\log(R/r)}.
$$
A little algebra shows that this is equal to
$$
\left(\frac{r}{R}\right)^{c\log(\alpha^{1-k/12})},
$$
and this is a bound of the required form, with $\lambda(k)= c\log(\alpha^{1-k/12})$.
\hfill $\Box$

\medskip\noindent
We now describe how the existence of subsequential weak limits of the sequence $\mu_{n}$ follows from this lemma. 
(This is well known but perhaps not immediately obvious from the literature, hence our summary for the convenience
of the reader.)  For $\epsilon >0$ and positive integers $k,n$ define

\begin{equation*}\label{def-crossing-rv}
r_{\epsilon,k}^{n}:=
 \inf \left( \left\{ 0<r \leq 1 : \begin{array}{c}
\text{ some annulus } A(x;r^{1+\epsilon},r), x \in [0,1]^2, \text{ is}\\
\text{crossed by } k \text{ disjoint crossings in } \FF_{n}
\end{array} \right\}, 1 \right) .
\end{equation*}

It follows exactly as in \cite{AizBur99} or \cite{Sun09} that, as a
consequence of Lemma \ref{lem-disjoint-crossings}, for any $\epsilon>0$, for large enough $k$ the random
variables $r_{\epsilon,k}^{n}$ are stochastically bounded away from
zero as $n \to \infty$, that is
\begin{equation}\label{bnd-stoch-away}
 \lim_{u \to 0}\limsup_{n \to \infty} \IP_p(r_{\epsilon,k}^{n} \leq u) =0.
\end{equation}
Note that in Lemma \ref{lem-disjoint-crossings} we have the result \eqref{bnd-disjoint-crossings}
only for $r$ small enough, while the corresponding hypothesis H1 in \cite{AizBur99}
is stated without that restriction. Our Lemma \ref{lem-disjoint-crossings}
however is sufficient to prove (\ref{bnd-stoch-away}), since the latter concerns the behaviour of $\IP_p(r_{\epsilon,k}^{n} \leq u)$
as $u \to 0$.

Note also that Hypothesis H1 is used in equation (3.4) in \cite{AizBur99} only for
annuli whose inner radius equals $3r_{n}^{1+\epsilon}$, with $r_n=2^{-n} \leq u$ for $u>0$.
By taking $u$ sufficiently small in the proof of Lemma 3.1 in \cite{AizBur99}, also
the inner radius $3r_n^{1+\epsilon}$ becomes sufficiently small and we can apply our
Lemma \ref{lem-disjoint-crossings}. Hence, (\ref{bnd-stoch-away}) follows from our
Lemma \ref{lem-disjoint-crossings} by the same arguments used in \cite{AizBur99} or \cite{Sun09}.

As shown in \cite{AizBur99} (see in particular their proof
of Theorem 1.1 and equation (1.7) in their Remark (ii) following the
statement of Theorem 1.1, but note the typo in equation (1.7), where
$d/\lambda(1)$ should be $\lambda(1)/d$), equation
(\ref{bnd-stoch-away}) implies the following result.

\begin{theorem}[\cite{AizBur99}]\label{thm-regularity}
For any $\epsilon > 0$, all curves $\Gamma \in \FF_{n}$
can be parametrized by continuous functions $\gamma: [0,1] \to [0,1]^2$
such that for each curve, for all $0 \leq t_1 \leq t_2 \leq 1$
$$
|\gamma(t_1)-\gamma(t_2)| \leq k_{\epsilon}^{n} |t_1 -
t_2|^{\frac{1}{2+\epsilon}},
$$
where the random variables $k_{\epsilon}^n$ are stochastically bounded
as $n \to \infty$, that is,
$$
\lim_{u \to \infty} \limsup_{n \to \infty} \IP_p(k^n_{\epsilon} \geq u) =0.
$$
\end{theorem}

Once we have this result, we use Theorem 5.7 in \cite{Sun09} and Theorem 1.1 in \cite{AizBur99} to conclude that
the sequence of measures $\{\mu_{n}\}_{n \geq 1}$ is tight. 
Since $\Sigma$ is separable, it then follows from Prohorov's theorem
that for every sequence $n_k \to \infty$ there exists a subsequence
$n_{k_{l}} \to \infty$ such that $\mu_{n_{k_{l}}}$ converges weakly to a
probability measure on $\Sigma$.

Finally, from the fact that the $k_{\epsilon}^n$ are stochastically bounded
and the fact that the collection of curves with a given H\"older exponent
is compact, we have (see also \cite{Sun09} and \cite{AizBur99}) that if
we sample from any such weak limit, all curves $\gamma$ in the sample can
be parametrized in such a way that
\begin{equation}
\label{eqnHolder}
|\gamma(t_1)-\gamma(t_2)| \leq M|t_1-t_2|^{\alpha},
\end{equation}
where $M$ is a random number common to all curves in the same sample,
and $\alpha$ is a (non-random) constant. 

Next we combine the above weak convergence with the a.s.\
convergence of the retained squares (as compact sets) in the fractal process.
Let $(S,H)$ denote the metric space of compact subsets of $[0,1]^d$ with
the Hausdorff distance $H$ and let $(\Sigma,{\rm Dist})$ be
as defined after equation (\ref{hausdorff-D}). Furthermore let the
function $g:\Sigma \mapsto S$ be defined by
$g({\cal F})=\cup_{\gamma \in {\cal F}} {\rm Image}(\gamma)$ and
define $F_n=g({\cal F}_n)$.

Our next result concerns weak convergence of
$({\cal F}_n,F_n)$ where we use the product topology on $\Sigma \times S$.

\begin{lemma} \label{lem-joint-weak-conv}
The distribution of $({\cal F}_n,F_n)$ converges weakly along a subsequence.
Furthermore, any pair $({\cal F},F)$ of random variables sampled from any
such weak limit a.s.\ satisfies $F=g({\cal F})$.
\end{lemma}

\medskip\noindent
{\em Proof.} We already know that the distribution $\mu_n$ of ${\cal F}_n$
converges weakly along a subsequence. Using (\ref{Distance})
and (\ref{hausdorff-D}), if ${\rm Dist}({\cal F},{\cal F}')\leq
\delta$, this immediately implies that $H(F,F')\leq \delta$ since
the images of any two curves $\gamma_1,\gamma_2$ such that
$D(\gamma_1,\gamma_2)\leq \delta$ are within Hausdorff distance
$\delta$. This proves that $g$ is continuous.

The convergence in distribution of ${\cal F}_n$ along some subsequence
$n_k$ to a limit $\cal F$ implies the existence of coupled versions $X_k$
and $X$ of ${\cal F}_{n_k}$ and $\cal F$ respectively, such that $X_k$ converges to
$X$ in probability as $k \to \infty$ (see, e.g. Corollary 1 in
\cite{billingsley71}). Moreover, since $g$ is continuous,
$g(X_k) \stackrel{dist.}{=} F_{n_k}$ converges in probability to
$g(X) \stackrel{dist.}{=} g({\cal F})$ as $k \to \infty$.
This implies convergence in probability of the vector $(X_k,g(X_k))$
to $(X,g(X))$, which yields the joint convergence in distribution of
$({\cal F}_{n_k},F_{n_k})$ to some limit $({\cal F},F)$ with $F=g({\cal F})$
a.s. \hfill $\Box$

\medskip\noindent
{\em Proof of Theorem \ref{thmnew1}.}
According to Lemma \ref{lem-joint-weak-conv} there exists a subsequence $\{n_k\}_{k \geq 1}$
such that $({\cal F}_{n_{k}},F_{n_{k}})$ converges weakly to some limit $({\cal F},F)$
where $F$ is a.s.\ the union of the images of ${\cal F}$. Furthermore,
we claim that a.s.
\begin{equation}
\label{as-convergence}
H(F_n,\C)\rightarrow 0 \, \text{ as } n \to \infty ,
\end{equation}
where $F_n=g(\FF_n)$. To see this, let $F_n^{\epsilon}$ denote the $\epsilon$-neighbourhood of
$F_n$ and note that $\C \not \subset F_n^{\epsilon}$ implies that there exists an $x \in \C$
such that $B(x,\epsilon) \cap F_n = \emptyset$ and so $B(x,\epsilon) \subset \C^n$,
otherwise an interface curve would be closer to $x$ than $\epsilon$. It is however easy to
prove that the probability that there exists a ball of radius $\epsilon$ in $\C^n$ goes to 0.
For the other direction, let $D=D(\C,\epsilon)$ be the complement of the open $\epsilon$-neighbourhood
of $\C$. Obviously $D \cap \C^n$ is a compact set (if nonempty) and furthermore
$D \cap \C_{n+1} \subset D \cap \C^n$ for every $n$. Therefore, by compactness, if
$D \cap \C^n \neq \emptyset$ for every $n$ then $\cap_{n=1}^{\infty} D \cap \C^n \neq \emptyset$
and so there are points in $\C$ that are also in $D$ which is a contradiction. Therefore, the
open $\epsilon$-neighbourhood of $\C$ will eventually contain $\C^n$ and hence also $F_n$.

Since $F_{n_k}$ converges weakly to $F$ and a.s.\ to $\C$ (because of (\ref{as-convergence})), we conclude that $F$ and
$\C$ have the same distribution, and hence $\C$ has the same distribution as $g(\FF)$. The fact that (as noted above) a.s., 
a realisation from $\FF$ only contains H\"older continuous curves satisfying (\ref{eqnHolder}) finishes the proof.
\hfill $\Box$

\medskip\noindent
{\em Proof of Theorem \ref{thm:lowest-crossing}.}
Let $\nu_n$ be the distribution of $\sigma_n$. 
We repeat the proof of Theorem \ref{thmnew1},
with the distance (\ref{hausdorff-D}) replaced by (\ref{Distance}). This shows that 
conditioned on the existence of a
left-right crossing for all $n$, $\nu_n$ has subsequential weak
limits as $n \to \infty$ and in addition that any such limit assigns probability 1 to H\"older
continuous curves. 
\hfill $\Box$

\section{Proof of Theorem \ref{thmnew2}}
\label{sec-Hausdorff}

We will start by showing that the fractal percolation process satisfies
Hypothesis H2 in \cite{AizBur99}, from which Theorem \ref{thmnew2} follows.
The hypothesis concerns probabilistic bounds on crossings in the long direction of
certain rectangles.

A collection of sets $\{A_i\}$ is called \emph{well-separated} if for all $i$ the distance of each set $A_i$ to
the other sets $\{A_j\}_{j\neq i}$ is at least as large as the diameter of $A_i$.
Hypothesis H2 in \cite{AizBur99} reads as follows.

\begin{hypothesis}\label{hyp-2} There exist $\sigma > 0$ and some $\rho<1$ such that
for every collection of $k$ well-separated rectangles $A_1,\ldots,A_k$ of width $l_1,\ldots,l_k$
and length $\sigma l_1,\ldots, \sigma l_k$, the following inequality holds:
 \begin{equation} \nonumber
  \limsup_{n \to \infty} \IP_p\left( \begin{array}{c}
\FF_{n} \text{ contains a long crossing in each of }
A_1,\ldots, A_k \end{array}
\right) \leq \rho^{k}.
\end{equation}
\end{hypothesis}

\begin{lemma}\label{lem-hyp-3}
Hypothesis \ref{hyp-2} holds for interface curves in the fractal process.
\end{lemma}

\medskip\noindent
{\em Proof.}
We assume without loss of generality that $l_1 \geq l_2 \geq \cdots \geq l_k$. Let $n_i$ be the smallest integer
$n$ for which {\em all} rectangles of dimensions $l_i$ and $\sigma l_i$ contain an $n$-level square of
the fractal construction, and define $\IA_i$ as the event of complete retention until iteration step $n_i$.
Let $CR^{I}_n(A_i)$ denote the event that the closed rectangle $A_i$ is crossed in the long direction
by an interface segment after $n$ iterations of the fractal process, and define $CR^{B}_n(A_i)$ similarly
for a black crossing.

Note that by the fact that $A_1,\ldots, A_k$ are well-separated and by the choice of $n_1$,
given $\IA_1$, the event $CR^{B}_n(A_1)$ is, when $n > n_1$, conditionally independent of
the events $CR^{B}_n(A_2),\ldots,CR^{B}_n(A_k)$, since no level-$n$ square intersects
more than one rectangle. Note also that, if a rectangle is crossed in the long direction by
an interface segment, then it also contains a black crossing in the long direction.
Using these two facts and the FKG inequality gives, for $n > n_1$,
\begin{eqnarray*}
 \IP_{p}(\cap_{i=1}^k CR^{I}_n(A_i)) &\leq& \IP_{p}(\cap_{i=1}^k CR^{B}_n(A_i)) \\
&\leq& \IP_{p}(\cap_{i=1}^k CR^{B}_n(A_i)  | \IA_{1}) \\
&=& \IP_{p}(CR^{B}_n(A_1)| \IA_1) \IP_{p}(\cap_{i=2}^k CR^{B}_n(A_i) | \IA_1).
\end{eqnarray*}
Since $l_1 \geq l_2 \geq \cdots \geq l_k$, we have $n_1 \leq n_2 \leq \cdots \leq n_k$ and hence
$\IA_1 \supset \IA_1 \supset \cdots \supset \IA_k$. It follows as before that
for $n > n_2$ we have
\begin{eqnarray*}
\IP_{p}(\cap_{i=2} ^k CR^{B}_n(A_i) | \IA_1)
&\leq& \IP_{p}(\cap_{i=2}^k CR^{B}_n(A_i) | \IA_2) \\
&=& \IP_{p}(CR^{B}_n(A_2) | \IA_2) \IP_{p}(\cap_{i=3}^k CR^{B}_n(A_i) | \IA_2).
\end{eqnarray*}
We repeat this procedure $k-2$ more times, finally obtaining, for $n > n_k$, that
\begin{equation} \nonumber
\IP_{p}(\cap_{i=1}^k CR^{I}_n(A_i)) \leq \prod_{i=1}^k \IP_{p}(CR^{B}_n(A_i)| \IA_i).
\label{eqn-conditioning-H2}
\end{equation}

It remains to show that $\IP_{p}(CR^{B}_n(A_i)| \IA_i)$ is uniformly bounded above by some $\rho < 1$.
This can be seen as follows. Let, for each $i$, $W_i$ be a smallest collection of
level-$(n_i+1)$ squares in the fractal process with the property that if each of the squares in $W_i$
is white, then $A_i$ cannot be crossed by a black path in the long direction. It is easy to see from
the choice of $n_i$ that the cardinality of $W_i$ is uniformly bounded above and hence that the
probability that all squares in $W_i$ are white is uniformly bounded below by some positive number.
It follows that the probability of a long black crossing in $A_i$ is uniformly bounded above by some
number strictly smaller than 1. \hfill $\Box$

\medskip\noindent
{\em Proof of Theorem \ref{thmnew2}.}
The result follows from Lemma \ref{lem-hyp-3} and Theorem 1.3 in \cite{AizBur99}. \hfill $\Box$

\medskip\noindent
\noindent{\bf Acknowledgments.}\, The second author thanks Lincoln
Chayes for suggesting the argument that led to the proof of
equation (\ref{eqn-main3.2}).

\end{document}